\documentclass[12pt]{amsart}
\usepackage{amsmath}
\usepackage{amsfonts}
\usepackage{amssymb}
\usepackage{graphicx}%
\providecommand{\U}[1]{\protect\rule{.1in}{.1in}}

\newtheorem{rmk}{Remark}

\newtheorem{prf}{Proof}

\newcommand{\NOT}[1]{}
\newcommand{\pa}{\par\medskip}

\title{The Role of Euclid's Circle of Ideas}

\author{Eliahu Levy}
\address{Department of Mathematics,
Technion -- Israel Institute of Technology,
Haifa 32000, Israel}
\email{eliahu@math.technion.ac.il}

\date{}

\keywords{Euclid, geometry, mathematics in Ancient Greece, axioms, experimental and observational evidence, Archimedes, Newton, Einstein, Hilbert, Kant, Born, The Axiom of Parallels, hyperbolic geometry, Cartesian coordinate system.}

\begin{document}
\maketitle


\begin{abstract}
In this small note I try to summarize some observations about Euclid's remarkable role in mathematics and about the ambient philosophy.
\end{abstract}

\medskip

The circle of ideas of Euclid's `Elements', the `hub' of the way to do mathematics that we associate with Ancient Greece, had an enormous, sometimes `dazzling' influence on the mathematical thought that followed. This dazzling is, of course, fully justified. Euclid has shown how a part of `reality' (indeed idealized -- but one may mention other idealizations in `proper' physics such as an `ideal gas'), is a mathematical system in the sense that the whole theory can be founded on a narrow basis of things whose knowledge is not a result of proofs/reasoning, on which a much larger edifice is built through proofs.\pa

In particular, such nontrivial facts as the three heights of a triangle intersecting in one point are derived from `pure thoughts and intuition' -- {\em a priori}, in the language of philosophers like Kant.\pa

Indeed, although any culture with mathematics developed to a stage where the usual methods to do arithmetics with numbers bigger than $10$, or the formula for the area of a rectangle, are recognized, would most likely justify them by reasoning, if it bothers to justify them at all, Euclid's tradition injected into mathematics the imperative to write a consciously explicitly stated proof, first in `geometry' and then in all mathematics.\pa

Later other parts of `reality' were found to also have such mathematical character --
you may deduce the whole edifice by `logic' from a narrow set of axioms --
Archimedes with levers and liquids, and after almost two millennia Newton with mechanics and more.\pa

But, historically, Euclid's geometry stood out. It seemed that here there is really `pure' mathematics, a `universe of discourse' of ideal points, lines etc.\ that `touches', but should be differentiated from `physics' -- `reality'. geometry should  associate itself primarily with other parts of mathematics -- numbers, magnitudes.\pa

As I said, people were dazzled. Kant, for one, made this a cornerstone of his efforts in philosophy.\pa

Of course, this had to interact with the developments and achievements of mathematics during the centuries; In the 17'th century Descartes' and Fermat's coordinates and analytic geometry; In the 19'th the non-Euclidean geometries and mathematics speaking about `possible universes', each with its own primitive notions and axioms.\pa

The physicist Max Born is quoted as saying that {\em physicists (and others) abhor axioms}. If one speaks about `reality' than an `axiom', in the sense of `something that is self evident, thus needs not any proof' is highly suspicious. If someone bases a logical edifice on such `axioms', they had better be derived from sound empirical evidence. Or at least, as with Newton's laws and others, stand the test of time -- i.e.\ the test of experiments and observations.\pa

But Euclid is not presented and conceived in that way. Thus Kant was convinced that it is `{\em a priori}' -- not based on empirical evidence at all.\pa

Let us have a closer look at Euclid and the `Elements'.\pa

The books of the `Elements' deal primarily with Plane Geometry, Numbers (arithmetics, with beautiful theorems like the existence of an infinite number of primes or the Euclidean algorithm to find the greatest common devisor of two natural numbers) and Space Geometry.\pa

Each subject begins with `definitions'. But only for Plane Geometry, in book I, there are `postulates' and `common notions' that traditionally are called also `axioms'. But these are, with two exceptions, rather general principles such as: `if one adds the same thing to equals one gets equals', or instructions to construct, such as `construct a circle with given center and radius' that one may interpret as a (constructive) statement of existence. For Numbers and Space Geometry there are only `definitions'.\pa

But what is highly striking is that while for Numbers and Space Geometry everything apparently follows just from the definitions, \textit{that is the case also with the sections about Plane Geometry}. With one exception, the `postulates' and `common notions' are only marginally used. In fact the main vehicle of the proofs are movements of figures, making them coinciding and making segments or triangles congruent, or one falling short of another etc.\pa

When Hilbert, around the end of the 19'th century, published his fully axiomatic version of Euclid (in his book `Foundations of Geometry'), he had to state basic congruence `theorems' as axioms (say, two triangles with equal corresponding sides and equal angles between them are congruent, thus equal in the other sides and angles) which {\em Euclid does not do}. For him, these are not axioms but consequences of the `moving' of figures.\pa

But the possibility of these movements is, albeit in an `idealized' form (in the sense of things like `ideal gas') a highly {\em empirical fact} -- in our world we can translate and rotate solid bodies -- which would not be possible if we lived in, say, the two dimensional surface of a sculpture (or if our space was curved, as general relativity teaches us to be the case, but in a negligible amount for all our usual purposes, from building houses to space travel).\pa

Euclid's geometry thus assumes a role, along with things like Newton's mechanics, of a valid part of physics, with assumptions (`axioms') justified empirically (either initially or by standing `the test of time').\pa

In particular (see Albert Einstein's famous `Geometrie und Erfahrung' article \cite{Einstein})
while for our everyday `rigid objects' we observe empirically the possibility of the
above-mentioned movements, to a good enough precision, thus `initially justifying Euclid', physics uses Euclidean geometry also in, say, atomic scales -- here the justification would not be initial but `standing the test of time'.\pa

Yet one must note the traditional difference of attitude. Firstly, contrary to Newton, Euclid succeeds in `concealing' the fact that the possibility of his `movements' rests on our empirical knowledge, thus `deceives' people like Kant. Secondly, contrary to mechanics etc.\ as a branch of physics, geometry, even in the `Elements', has its eye on other mathematics, numbers etc, hence mathematical rigor is always `at the horizon'. A book like Hilbert's would be much less expected in mechanics.\pa

The real reason for that is that geometry, even as an algebra of coordinates, has a rich `intrinsic' mathematical contents. When mechanics has also something to contribute here, such as the notion of center of equilibrium or densities of `weight', they fall naturally under forms of `geometry', inasmuch as one is concerned with `mathematics' intrinsic contents'.\pa

But let us return to the history. The great `crisis' of a sort concerning Euclid came
from the one real axiom in his book -- the fifth postulate there -- the axiom of parallels. (There is another `real axiom' -- the fourth postulate, stating that all right angles are equal, but it can easily be proved to Euclid's school's satisfaction.) The fifth postulate is an {\em axiom to be abhorred}. It is clear that it is not `self evident' in any way. On the other hand, many theorems depend on it, to mention just the fact that the sum of angles in a triangle is $180^\circ$. Also, it conforms to experiment in our mundane world (we know now that it will not when General Relativity becomes significant).\pa

Thus, for two thousand years people tried to prove it, i.e.\ remove it from the status of `abhorrent axiom'. Then in the beginning of the 19'th century Gauss, Bolyai and Lobachevsky had the `revolutionary' idea that maybe it cannot be proved -- one may proceed to construct a geometry where it is incorrect. This turned to be the hyperbolic geometry, one of the non-Euclidean geometries -- a mathematical structure of great importance quite besides this `axiomatic' issue. It is also homogeneous -- one may make movements freely (another homogenous geometry is the geometry on the sphere, which may also be considered a non-Euclidean geometry). Indeed, the parts of the `Elements' that rely only on the `movements' and do not use the axiom of parallels are fully valid also in the hyperbolic geometry.\pa

In our times, a somewhat similar attitude/role to that historically  of Euclid is enjoyed by basing our geometry on the Cartesian coordinate system, or, more neatly expressed, using vectors etc.\ --  for `modern' science, engineering etc.\ that is the `Euclidean geometry'. In this way, Euclidean geometry is a branch of algebra, of course, with no special role to an Axiom of Parallels. But the conviction that this is indeed our physical space -- real things will `obey' what is obtained mathematically in this `structure', lies on a similar basis as the old conviction in the classical Euclid.

\end{document}